# Proper Affine Vector Fields in Plane Symmetric Static Space-Times

Ghulam Shabbir

Faculty of Engineering Sciences

GIK Institute of Engineering Sciences and Technology

Topi Swabi, NWFP, Pakistan

Email: shabbir@giki.edu.pk

and

Nisar Ahmed

Faculty of Electronic Engineering

GIK Institute of Engineering Sciences and Technology

Topi Swabi, NWFP, Pakistan.

**Abstract:**

A study of proper affine vector fields in plane symmetric static space-times by using the rank of the $6\times 6$ Rieman matrix and holonomy. Studying proper affine vector fields in each case, It is shown that the special class of the above space-times admit proper affine vector fields.

**1. INTRODUCTION**

Over the recent years, there has been much interest in the study of symmetry in general relativity; affine symmetry seems to ignored probably because the complications involved in (1). Following [1,9] some progresss was made in [2] where they actually shown that under what conditions the space-times admit proper affine vector fields. Other different approaches may be seen in [9,10,11]. In this paper a different approach is devolped to study proper affine vector fields in plane symmetric static space-times by using holonomy, the rank of the $6\times 6$ Rieman matrix and direct integration techniques. Through out M is representing the four dimensional, connected, hausdorff space-time manifold with Lorentz metric g of signature (-, +, +, +). The curvature tensor associated with g, through Levi-Civita connection, is denoted in component form by $R^{a}{}_{bcd}$. The usual covariant, partial and Lie derivatives are denoted by a semicolon, a comma and the symbol L, respectively. Round and square brackets denote the usual symmetrization and skew-symmetrization, respectively. The space-time $M$ will be assumed nonflat in the sense that the Riemann tensor does not vanish over any nonempty open subset of $M$.

A vector field $X$ on $M$ is called an affine vector field if it satisfies

$$X_{a;bc} = R_{abcd} X^d. \tag{1}$$

If one decomposes $X_{a;b}$ on $M$ into its symmetric and skew-symmetric parts

$$X_{a;b} = \frac{1}{2} h_{ab} + F_{ab} \qquad (h_{ab} = h_{ba}, \quad F_{ab} = -F_{ba}) \tag{2}$$

then equation (1) is equivalent to

$$(i)\ h_{ab;c} = 0 \quad (ii)\ F_{ab;c} = R_{abcd} X^d \quad (iii)\ F_{ab;c} X^c = 0. \tag{3}$$

The above vector field $X$ is called an affine on vector field $M$ if the local diffeomorphisms $\psi_t$ (for appropriate $t$) associated with $X$ map geodesics into geodesics and also preserve their affine parameter. If $h_{ab} = 2c g_{ab}$, $c \in \Re$, then the vector field $X$ is called homothetic (and *Killing* if $c = 0$). The vector field $X$ is said to be proper affine if it is not homothetic vector field and also $X$ is said to be proper homothetic vector field if it is not Killing vector field on $M$ [1]. It is important to note the following equation

$$R_{abcd} k^d = 0. \tag{4}$$

## 2. Affine Vector Fields

In this section we will briefly discuss when the space-times admit proper affine vector fields for further details see [2,8].

Suppose that $M$ is a simple connected space-time. Then the holonomy group of $M$ is a connected Lie subgroup of the idenity component of the Lorentz group and is thus characterized by its subalgebra in the Lorentz algebra. These have been labeled into fifteen types $R_1 - R_{15}$ [3,4]. It follows from [2] that the only such space-times which could admit proper affine vector fields are those which admit nowhere zero covariantly constant second order symmetric tensor field $h_{ab}$ and it is known that this forces the holonomy type to be either $R_2$, $R_3$, $R_4$, $R_6$, $R_7$, $R_8$, $R_{10}$, $R_{11}$ or $R_{13}$. Here, we will only discuss the space-times which has the holonomy type $R_2$, $R_4$, $R_7$, $R_{10}$ or $R_{13}$.

First consider the case when $M$ has type $R_{13}$. Then one can always set up local coordinates $(t, x^1, x^2, x^3)$ on an open set $U = U_1 \times U_2$, where $U_1$ is a one dimensional timelike submanifold of $U$ coordinatized by $t$ and $U_2$ is a three dimensional spacelike submanifold of $U$ coordinatized by $x^1, x^2, x^3$ and where the above product is a metric product and the metric on $U$ is given by [1]

$$ds^2 = -dt^2 + g_{\alpha\beta} dx^\alpha dx^\beta \qquad (\alpha, \beta = 1, 2, 3) \qquad (5)$$

where $g_{\alpha\beta}$ depends on $x^\gamma$, ($\gamma = 1, 2, 3$). The above space-time is clearly 1+3 decomposable. The curvature rank of the above space-time is atmost three and there exists a unique nowhere zero vector field $t_a = t_{,a}$ satisfying $t_{a;b} = 0$ and also $t_a t^a = -1$. From the Ricci Identity $R^a{}_{bcd} t^d = 0$. It follows from [2] that affine vector fields in this case are

$$X = (c_1 t + c_2) \frac{\partial}{\partial t} + Y \qquad (6)$$

where $c_1, c_2 \in \Re$ and $Y$ is a homothetic vector field in the induced geometry on each of the three dimensional submanifolds of constant $t$.

Now consider the situation when $M$ has type $R_{10}$. The situation is similar to that of previous $R_{13}$ case except that now we have local decomposition is $U = U_1 \times U_2$, where $U_1$ is a one dimensional spacelike submanifold of $U$ and $U_2$ is a three dimensional timelike submanifold of $U$. The space-time metric on $U$ is given by [1]

$$ds^2 = dx^2 + g_{\alpha\beta} dx^\alpha dx^\beta \qquad (\alpha, \beta = 0, 2, 3) \qquad (7)$$

where $g_{\alpha\beta}$ depends on $x^\gamma$, ($\gamma = 0, 2, 3$). The above space-time is clearly 1+3 decomposable. The curvature rank of the above space-time is atmost three and there exists a unique nowhere zero vector field $x_a = x_{,a}$ satisfying $x_{a;b} = 0$ and also $x_a x^a = 1$. From the Ricci Identity $R^a{}_{bcd} x^d = 0$. It follows from [2] that affine vector fields in this case are

$$X = (c_1 x + c_2)\frac{\partial}{\partial x} + Y \tag{8}$$

where $c_1, c_2 \in \Re$ and $Y$ is a homothetic vector field in the induced geometry on each of the three dimensional submanifolds of constant $x$.

Next suppose $M$ has type $R_7$. Then each $p \in M$ has a neighborhood $U$ which decomposes metrically as $U = U_1 \times U_2$, where $U_1$ is a two dimensional submanifold of $U$ with an induced metric of Lorentz signature and $U_2$ is a two dimensional submanifold of $U$ with positive definite induced metric. The space-time metric on $U$ is given by [1]

$$ds^2 = P_{AB} dx^A dx^B + Q_{\alpha\beta} dx^\alpha dx^\beta \tag{9}$$

where $P_{AB} = P_{AB}(x^C), \forall A,B,C = 0,1$ and $Q_{\alpha\beta} = Q_{\alpha\beta}(x^\gamma), \forall \alpha,\beta,\gamma = 2,3$ and the above space-time is clearly 2+2 decomposable. The space-time (8) admits two recurrent vector fields [5] $l$ and n i.e. $l_{a;b} = l_a p_b$ and $n_{a;b} = n_a p_b$ where $p_b$ is the recurrent 1-form. It also admits two covariantly constant second order symmetric tensors which are $2l_{(a} n_{b)}$ and $(x_a x_b + y_a y_b)$. The rank of the $6 \times 6$ Riemann matrix is two. It follows from [2] that if $X$ is an affine vector field on $M$ then $X$ decomposes as

$$X = X_1 + X_2 \tag{10}$$

where the vector fields $X_1$ and $X_2$ are tangent to the two dimensional timelike and spacelike submanifolds, respectively. It also follows from [2] that $X_1$ and $X_2$ are homothetic vector fields in their respective submanifolds with their induced geometry. Conversely, every pair of affine vector fields, one in the timelike submanifolds and one spacelike submanifolds give rise to a affine vector field in space-time.

Now suppose that $M$ has type $R_4$. Then each $p \in M$ has a neighborhood $U$ which decomposes metrically as $U = U_1 \times U_2 \times U_3$, where $U_1$ and $U_2$ are one dimensional submanifold of $U$ and $U_3$ is a two dimensional submanifold of $U$. The space-time metric on $U$ is given by [2]

$$ds^2 = -dt^2 + dx^2 + g_{AB} dx^A dx^B \qquad (A, B = 2,3) \qquad (11)$$

where $g_{AB}$ depends only on $x^C$ ($C = 2, 3$). The above space-time is clearly 1+1+2 decomposable. The curvature rank of the above space-time is one and there exist two independent nowhere zero unit timelike and spacelike covariantly constant vector field $t_a = t_{,a}$ and $x_a = x_{,a}$ satisfying $t_{a;b} = 0$ and $x_{a;b} = 0$. From the Ricci identity $R^a{}_{bcd} t_a = 0$ and $R^a{}_{bcd} x_a = 0$. It follows from [2] that affine vector fields in this case are

$$X = (c_1 t + c_2 x + c_3) \frac{\partial}{\partial t} + (c_4 t + c_5 x + c_6) \frac{\partial}{\partial x} + Y \qquad (12)$$

where $c_1, c_2, c_3, c_4, c_5, c_6 \in \Re$ and $Y$ is a homothetic vector on each of two dimensional submanifolds of constant $t$ and $x$.

Now suppose that $M$ has type $R_2$. Here each $p \in M$ admits a neighborhood $U$ which decomposes metrically as $U = U_1 \times U_2 \times U_3$, where $U_1$ and $U_2$ are one dimensional submanifold of $U$ and $U_3$ is a two dimensional submanifold of $U$. The space-time metric on $U$ is given by [2]

$$ds^2 = dy^2 + dz^2 + g_{AB} dx^A dx^B \qquad (A, B = 0,1) \qquad (13)$$

where $g_{AB}$ depends only on $x^C$ ($C = 0,1$). The above space-time is clearly 1+1+2 decomposable. The curvature rank of the above space-time is one and there exist two independent nowhere zero unit timelike and spacelike covariantly constant vector field $y_a = y_{,a}$ and $z_a = z_{,a}$ satisfying $y_{a;b} = 0$ and $z_{a;b} = 0$. From the Ricci identity $R^a{}_{bcd} y_a = 0$ and $R^a{}_{bcd} z_a = 0$. It follows from [2] that affine vector fields in this case are

$$X = (c_1 y + c_2 z + c_3) \frac{\partial}{\partial y} + (c_4 y + c_5 z + c_6) \frac{\partial}{\partial z} + Y \qquad (14)$$

where $c_1, c_2, c_3, c_4, c_5, c_6 \in \Re$ and $Y$ is a homothetic vector on each of two dimensional submanifolds of constant $y$ and $z$.

## 3. MAIN RESULTS

As mentioned in section 2, the space-times which can admit proper affine vector fields having holonomy type $R_2$, $R_3$, $R_4$, $R_6$, $R_7$, $R_8$, $R_{10}$, $R_{11}$ or $R_{13}$. It also follows from [5] that the rank of the $6 \times 6$ Riemann matrix is atmost three. Here in this paper we will consider the rank of the $6 \times 6$ Riemann matrix to study affine vector fields in plane symmetric static space-time. Consider plane symmetric static space-time in the usual coordinate system $(t, x, y, z)$ with line elememt [7]

$$ds^2 = -e^{\nu(x)} dt^2 + dx^2 + e^{\mu(x)}(dy^2 + dz^2) \tag{15}$$

It follows from [7], the above space-time admits four independent Killing vector fields which are

$$\frac{\partial}{\partial t}, \quad \frac{\partial}{\partial y}, \quad y\frac{\partial}{\partial z} - z\frac{\partial}{\partial y}, \quad \frac{\partial}{\partial z}.$$

The non-zero independent components of the Riemann tensor are

$$\begin{aligned}
R^{01}{}_{01} &= \frac{1}{4}(2\nu'' + \nu'^2)e^{-\nu} \equiv \alpha_1, \\
R^{02}{}_{02} &= R^{03}{}_{03} = -\frac{1}{4}\mu'\nu' \equiv \alpha_2, \\
R^{12}{}_{12} &= R^{13}{}_{13} = -\frac{1}{4}(2\mu'' + \mu'^2) \equiv \alpha_3, \\
R^{23}{}_{23} &= -\frac{1}{4}\mu'^2 \equiv \alpha_4.
\end{aligned} \tag{16}$$

Writing the curvature tensor with components $R^{ab}{}_{cd}$ at p as a 6X6 symmetric matrix in a well known way [6]

$$R^{ab}{}_{cd} = diag(\alpha_1, \alpha_2, \alpha_2, \alpha_3, \alpha_3, \alpha_4)$$

where $\alpha_1, \alpha_2, \alpha_3$ and $\alpha_4$ are real functions of $x$. The 6-dimensional labeling is in the order 01, 02, 03, 12, 13, 23 with $x^0 = t$. We are only interested in those case when the rank of the $6 \times 6$ Riemann matrix is less than or equal to three (excluding flat cases). We thus obtain the following cases:

(A) Rank $\leq 3$, when $\nu \in R$, $\mu = \mu(r)$,

(B) Rank $= 1$, when $\nu = \nu(x), \lambda \in R$,

(C) Rank =3, when $v = v(x)$, $\mu = \mu(x)$, $2v'' + v'^2 = 2\mu'' + \mu'^2 = 0$.

We consider each case in turn.

## Case A

In this case $v \in R, \mu = \mu(r)$. First suppose that the rank of the $6 \times 6$ Riemann matrix is 3 and there exists a unique (up to a multiple) nowhere zero timelike vector field $t_a = t_{,a}$ such that $t_{a;b} = 0$ (and so, from the Ricci idenity $R^a{}_{bcd} t_a = 0$). The line element can, after a recaling of t, be written in the form

$$ds^2 = -dt^2 + (dx^2 + e^{\mu(x)}(dy^2 + dz^2)). \qquad (17)$$

The space-time is clearly 1+3 decomposable. The affine vector fields in this case [2] are

$$X = (c_8 t + c_9)\frac{\partial}{\partial t} + X' \qquad (18)$$

where $c_8, c_9 \in R$ and $X'$ is a homothetic vector field in the induced geometry on each of the three dimensional submanifolds of constant $t$. The completion of case A necessities finding an homothetic vector fields in the induced geometry of the submanifolds of constant $t$. The induced metric $g_{\alpha\beta}$ (where $\alpha, \beta = 1, 2, 3$) with nonzero components is given by

$$g_{11} = 1, \ g_{22} = e^{\mu(x)}, \ g_{33} = e^{\mu(x)}. \qquad (19)$$

A vector field $X'$ is called homothetic vector field if it satisfies $L_{X'} g_{\alpha\beta} = 2c g_{\alpha\beta}$, where $c \in R$. One can expand by using (19) to get

$$X^1{}_{,1} = c \qquad (20)$$

$$X^1{}_{,2} + e^{\mu(x)} X^2{}_{,1} = 0 \qquad (21)$$

$$X^1{}_{,3} + e^{\mu(x)} X^3{}_{,1} = 0 \qquad (22)$$

$$\mu' X^1 + 2X^2{}_{,2} = 2c \qquad (23)$$

$$X^2{}_{,3} + X^3{}_{,2} = 0 \qquad (24)$$

$$\mu' X^1 + X^3{}_{,3} = 2c. \qquad (25)$$

From equations (20), (21) and (22) we get

$$X^1 = cx + A^1(y,z), \quad X^2 = -A^1{}_y(y,z)\int e^{-\mu}dx + A^2(y,z)$$

$$X^3 = -A^1{}_z(y,z)\int e^{-\mu}dx + A^3(y,z) \tag{26}$$

where $A^1(y,z)$ $A^2(y,z)$ and $A^3(y,z)$ are functions of integration. Now consider equation (24) and differentiate with respect to x and using (26) one finds $A^1{}_{yz}(y,z) = 0 \Rightarrow A^1(y,z) = B^1(y) + B^2(z)$ where $B^1(y)$ and $B^2(y)$ are functions of integration. Substituting back into (26) we get

$$X^1 = cx + B^1(y) + B^2(z), \quad X^2 = -B^1{}_y(y)\int e^{-\mu}dx + A^2(y,z),$$

$$X^3 = -B^2{}_z(z)\int e^{-\mu}dx + A^3(y,z). \tag{27}$$

Substituting equation (25) from (23) and differentiate with respect to x and using (27) gives $B^1{}_{yy}(y) - B^2{}_{zz}(z) = 0$, differentiate with respect to y gives $B^1{}_{yyy}(y) = 0 \Rightarrow B^1(y) = \frac{c_1}{2}y^2 + c_2 y + c_3$. Substituting back we get $B^2(z) = \frac{c_1}{2}z^2 + c_4 z + c_5$ where $c_1, c_2, c_3, c_4, c_5 \in \Re$ and so equation (27) becomes

$$X^1 = cx + (\frac{c_1}{2}y^2 + c_2 y + c_3) + (\frac{c_1}{2}z^2 + c_4 z + c_5),$$

$$X^2 = -(c_1 y + c_2)\int e^{-\mu}dx + A^2(y,z),$$

$$X^3 = -(c_1 z + c_4)\int e^{-\mu}dx + A^3(y,z). \tag{28}$$

Now consider equation (23), taking partial derivatives with respect to y and using (28) to get $\mu'(c_1 y + c_2) + 2A^2{}_{yy}(y,z) = 0$. Differentiating with respect to x gives $\mu''(c_1 y + c_2) = 0$ and there now exist two possibilities:

(1) $\mu'' = 0, \quad (c_1 y + c_2) \neq 0,$ (2) $\mu'' \neq 0, \quad (c_1 y + c_2) = 0.$

**Class A1**

$\mu'' = 0 \Rightarrow \mu = ax + b,$ where $a,b \in R (a \neq 0)$ and $(c_1 y + c_2) \neq 0$. Using this information in (28) and solving equations (20) to (25), one finds that homothetic vector fields in this case are Killing vector fields which are given by

$$X^1 = c_2 y + c_4 z + c_5,$$

$$X^2 = c_2 \frac{1}{a} e^{-ax-b} + c_2 \frac{a}{4}(z^2 - y^2) - \frac{a}{2}(c_3 yz + c_4 y) + c_5 z + c_6,$$

$$X^3 = c_2 \frac{1}{a} e^{-ax-b} + c_3 \frac{a}{4}(-z^2 + y^2) - \frac{a}{2}(c_2 yz + c_4 y) - c_5 z + c_7 \qquad (29)$$

where $c_2, c_3, c_4, c_5, c_6, c_7 \in R$. Thus the submanifolds of constant $t$ are of constant curvature and affine vector fields in this case are

$$X^0 = c_8 t + c_9, \qquad X^1 = c_2 y + c_4 z + c_5,$$

$$X^2 = c_2 \frac{1}{a} e^{-ax-b} + c_2 \frac{a}{4}(z^2 - y^2) - \frac{a}{2}(c_3 yz + c_4 y) + c_5 z + c_6,$$

$$X^3 = c_2 \frac{1}{a} e^{-ax-b} + c_3 \frac{a}{4}(-z^2 + y^2) - \frac{a}{2}(c_2 yz + c_4 z) - c_5 y + c_7. \qquad (30)$$

## Case A2

In this case $\mu'' \neq 0$ and $(c_1 y + c_2) = 0 \Rightarrow c_1 = c_2 = 0$ and $A^2{}_{yy}(y,z) = 0$ or $A^2(y,z) = yB^3(z) + B^4(z)$ where $B^3(z)$ and $B^4(z)$ are functions of integrations. Using the information in (28) and solving equations (20) to (25), one finds that there exist the following two possible subcasses:

(i) $\quad 2\mu'' + \mu'^2 \neq 0$ \qquad (ii) $\quad 2\mu'' + \mu'^2 = 0$

## Case A2i

In this case $2\mu'' + \mu'^2 \neq 0$ and that the rank of the $6 \times 6$ Riemann matrix is 3. Using the above information in equations (20) to (25) and solving them, one finds that homothetic vector fields in this case are Killing vector fields which are given by

$$X^1 = 0, \qquad X^2 = -c_1 z + c_3,$$

$$X^3 = c_1 y + c_2. \qquad (31)$$

where $c_1, c_2, c_3 \in R$ and affine vector fields in this case are (from equation (18))

$$X^0 = c_8 t + c_9, \; X^1 = 0, \qquad X^2 = -c_1 z + c_3,$$

$$X^3 = c_1 y + c_2.$$

## Case A2ii

In this case $2\mu'' + \mu'^2 = 0 \Rightarrow \mu = \ln(ax+b)^2$, where $a,b \in R\,(a \neq 0)$. Here the rank of the $6 \times 6$ Riemann matrix is one and there exist two independent nowhere zero solutions of equation (4) i.e. $R^a{}_{bcd} t_a = R^a{}_{bcd} x_a = 0$ where $x_a = x_{,a}$ and $t_a = t_{,a}$ are the spacelike and timelike vector field, respectively. The sapce-time (32) admits only one independent nowhere zero timelike covariantly constant vector field $t_a$ satisfying $t_{a;b} = 0$. After rescaling of $t$ the line element takes the form

$$ds^2 = dt^2 + dx^2 + (ax+b)^2(dy^2 + dz^2). \tag{32}$$

The above space-time is clearly 1+3 decomposable but the rank of 6X6 Riemann matrix is one. Substituting the above information into affine equations one finds, affine vector fields in this case

$$X^0 = c_4 t + c_5 x + c_6, \quad X^1 = c_7 t + c_8 x, \quad X^2 = -c_1 z + c_3,$$

$$X^3 = c_1 y + c_2 \tag{33}$$

where $c_1, c_2, c_3, c_4, c_5, c_6, c_7, c_8 \in R$. This completes case (A).

## Case B

In this case $\mu \in R$, $v = v(x)$ and the rank of the $6 \times 6$ Riemann matrix is one and there exist two independent nowhere zero spacelike vector field $y_a = y_{,a}$ and $z_a = z_{,a}$ satisfying $y_{a;b} = 0$ and $z_{a;b} = 0$. From the Ricci identity $R^a{}_{bcd} y_a = R^a{}_{bcd} z_a = 0$. The line element can, after a rescaling of $y$ and $z$, be written as

$$ds^2 = dy^2 + dz^2 + (-e^{v(x)} dt^2 + dx^2). \tag{34}$$

Clearly the above space-time is 1+1+2 decomposable. The affine vector fields in this case are [2]

$$X = (c_1 y + c_2 z + c_3)\frac{\partial}{\partial y} + (c_4 y + c_5 z + c_6)\frac{\partial}{\partial z} + X' \tag{35}$$

where $c_1, c_2, c_3, c_4, c_5, c_6 \in R$ and $X'$ is a homothetic vector fields in each of two dimensional submanifolds of constant $y$ and $z$. The next step is to finding homothetic vector fields in the induced geometry of the submanifolds of constant $y$ and $z$. The induced metric $g_{AB}$ (where $A, B = 0,1$) with nonzero components is given by

$$g_{00} = -e^{v(x)}, \ g_{11} = 1. \tag{36}$$

A vector field $X'$ is called homothetic vector field if it satisfies $L_{X'} g_{AB} = 2c g_{AB}$, where $c \in R$. One can expand by using (36) to get

$$v' X^1 + 2 X^0{}_{,0} = 2c \tag{37}$$

$$X^1{}_{,0} - e^v X^0{}_{,1} = 0 \tag{38}$$

$$X^1{}_{,1} = c. \tag{39}$$

Equation (39) gives $X^1 = cx + A^1(t)$, where $A^1(t)$ is a function of integration. Using value of $X^1$ in equation (38) we get $X^0 = A^1_t(t) \int e^{-v} dx + A^2(t)$, where $A^2(t)$ is a function of integration. If one proceeds further, after a straightforward calculation one finds that proper homothetic vector field exist if and only if $v = \ln(ax^2)$, where $a \in R - \{0\}$. Substituting the value of $v$ into (16), one finds that the rank of $6 \times 6$ Riemann matrix reduces to zero thus giving a contradiction (since we are assuming that the rank of $6 \times 6$ Riemann matrix is one). So homothetic vector fields in the induced geometry of constant $y$ and $z$ are Killing vector fields. If one proceeds further one finds there exist two possibilities:

(1) $v'' e^v = c$   (2) $v'' e^v \neq c$

where $c \in R$.

## Case B1

In this case further three possibilities exist

(i) $c > 0$,    (ii) $c < 0$,    (iii) $c = 0$.

We will consider each case in turn.

**(i)** Affine vector fields in this case are

$$X^0 = \sqrt{\frac{c}{2}}(c_7 \cos\sqrt{\frac{c}{2}}t - c_8 \sin\sqrt{\frac{c}{2}}t)\int e^{-v}dx + c_9,$$

$$X^1 = c_7 \sin\sqrt{\frac{c}{2}}t + c_8 \cos\sqrt{\frac{c}{2}}t \qquad (40)$$

$$X^2 = c_1 y + c_2 z + c_3, \qquad X^3 = c_4 y + c_5 z + c_6.$$

provided that $v''e^v = c$. Where $c_7, c_8, c_9 \in R$.

**(ii)** In this case $c < 0$. Put $c = -N$, where $N \in R (N > 0)$. Affine vector fields in this case are

$$X^0 = \sqrt{\frac{N}{2}}(c_7 \cosh\sqrt{\frac{N}{2}}t + c_8 \sinh\sqrt{\frac{N}{2}}t)\int e^{-v}dx + c_9,$$

$$X^1 = c_7 \sinh\sqrt{\frac{N}{2}}t + c_8 \cosh\sqrt{\frac{N}{2}}t \qquad (41)$$

$$X^2 = c_1 y + c_2 z + c_3, \qquad X^3 = c_4 y + c_5 z + c_6.$$

provided that $v''e^v = -N$. Where $c_7, c_8, c_9 \in R$.

**(iii)** In this case $c = 0 \Rightarrow v = ax + b$, where $a, b \in R (a \neq 0)$. Affine vector fields in this case are

$$X^0 = -\frac{c_7}{a} e^{-(ax+b)} + c_9 \qquad X^1 = c_7 t + c_8 \qquad (42)$$

$$X^2 = c_1 y + c_2 z + c_3, \qquad X^3 = c_4 y + c_5 z + c_6.$$

where $c_7, c_8, c_9 \in R$.

## Case B2

Affine vector fields in this case are

$$X^0 = c_7, \qquad X^1 = 0, \qquad (43)$$

$$X^2 = c_1 y + c_2 z + c_3, \qquad X^3 = c_4 y + c_5 z + c_6.$$

where $c_1 \in R$. This completes case B.

## Case C

In this case $v = v(r)$, $\mu = \mu(x)$, $2v'' + v'^2 = 0$ and $2\mu'' + \mu'^2 = 0$. The equations $2\mu'' + \mu'^2 = 0$ and $2v'' + v'^2 = 0 \Rightarrow v = \ln(\frac{a}{2}x + b)^2$ and $2\mu'' + \mu'^2 = 0 \Rightarrow \mu = \ln(\frac{c}{2}x+d)^2$ where $a,b,c,d \in R (a,c \neq 0)$ and $a \neq c, b \neq d$

(by assumption). The rank of the $6 \times 6$ Riemann matrix is 3 and there exist a unique (up to a multiple) nowhere zero spacelike vector field $x_a = x_{,a}$ such that $R^a{}_{bcd} x_a = 0$ and $x_{a;b} \neq 0$. The line element is

$$ds^2 = -(\frac{a}{2}x+b)^2 dt^2 + dx^2 + (\frac{c}{2}x+d)^2(dy^2 + dz^2) \qquad (44)$$

Substituting the above information into affine equations and after a strightforward calculation one find affine vector fields in this case are

$$X^0 = c_2, \qquad X^1 = 0, \qquad X^2 = c_3 z + c_4,$$
$$X^3 = -c_3 y + c_5 \qquad (45)$$

where $c_2, c_3, c_4, c_5 \in \Re$. Clearly affine vector fields in this case are Killing vector fields.

Now consider the case when $a = c \neq 0, b = d$, the line element is

$$ds^2 = -(\frac{a}{2}x+b)^2 dt^2 + dx^2 + (\frac{a}{2}x+b)^2(dy^2 + dz^2). \qquad (46)$$

It follows from [9] affine vector fields in this case are

$$X^0 = c_1 y + c_2 z + c_3, \quad X^1 = c_7 x,$$
$$X^2 = c_1 t - c_4 z + c_5, \quad X^3 = c_2 t + c_4 z + c_6 \qquad (47)$$

Where $c_1, c_2, c_3, c_4, c_5, c_6, c_7 \in R$.

## SUMMARY

In this paper a study of plane symmetric static space-times according to their proper affine vector fields is given. An approach is developed to study proper affine vector fields in the above space-times by using the rank of the 6X6 Riemann matrix, holonomy and direct integration techniques. From the above study we obtain the following results:

(i)   The case when the rank of the 6X6 Riemann matrix is three and there exists a nowhere zero independent spacelike vector field which is the solution of equation (4) and is not covariantly constant. This is the space-time (44) and it admits affine vector fields which are Killing vector fields (see for details Case C).

(ii)    The case when the rank of the 6X6 Riemann matrix is three and there exist a nowhere zero independent spacelike vector field which is a solution of equation (4) and is not covariantly constant. This is the space-time (46) and it admits the proper affine vector fields (equation 47).

(iii)   The case when the rank of the 6X6 Riemann matrix one there exist two nowhere zero independent solutions of equation (4) but only one independent nowhere zero covariantly constant vector field. This is the space-time (32) and it admits the proper affine vector fields (see Case A2ii).

(iv)    The case when the rank of the 6X6 Riemann matrix is two or three and there exists a nowhere zero independent timelike vector field which is the solution of equation (4) and also covariantly constant. This is the space-time (17) and it admits the proper affine vector fields (see Cases A1 and A2i).

(v)     The case when the rank of the 6X6 Riemann matrix one there exist two nowhere zero independent spacelike vector fields which are solutions of equation (4) and which are covariantly constant vector field. This is the space-time (34) and it admits the proper affine vector fields (see Case B1i, B1ii, B1iii and B2).